\documentclass[twoside,reqno,a4paper,12pt]{amsart}
\usepackage{amsmath,amsthm}
\usepackage{amsmath,amstext,amsfonts}
\usepackage{hyperref}
\usepackage{graphicx}
\usepackage{graphics}
\usepackage{tikz}
\usepackage{float}
\usepackage{etex}

\theoremstyle{plain}

\newtheorem{lem}{Lemma}

\theoremstyle{remark}

\DeclareMathOperator{\Real}{Re}

\newcommand{\R}{\mathbb R}

\newcommand{\cqd}{\hfill $\rule{2.5mm}{2.5mm}$}

\def \to {\rightarrow}

\def\cqd{\hfill{\vrule height 1.5ex width 1.4ex depth -.1ex}}

\newcommand{\trip}[1]{{|\kern -1pt\|#1\|\kern -1pt|}}

\begin{document}

\title[A new example on Lyapunov stability]
{A new example on Lyapunov stability}

\author[Hildebrando M. Rodrigues]
{Hildebrando M. Rodrigues$^{\dag}$}

\thanks{\dag Instituto de Ciências Matemáticas e de
Computaç\~{a}o, Universidade de S\~{a}o Paulo-Campus de S\~{a}o
Carlos, Caixa Postal 668, 13560-970 S\~{a}o Carlos SP, Brazil, e-mail: hmr@icmc.usp.br}

\author[J. Solà-Morales]
{J. Sol\`a-Morales$^{\ddag}$}
\thanks{ \ddag Departament de
Matemàtiques, Universitat Politècnica de Catalunya,
Av. Diagonal 647, 08028 Barcelona, Spain, e-mail:
jc.sola-morales@upc.edu}
 \markboth{Joan de Sol\`a-Morales}{
Hildebrando M. Rodrigues}
\thanks{\dag Partially supported by FAPESP Processo 2018/05218-8 and 
PQ-Sr 2018-Bolsa de Produtividade em pesquisa SÊNIOR
Processo CNPq 304767/2018-2 }
\thanks{\ddag Partially supported by MINECO grant MTM2017-84214-C2-1-P. Faculty member of the Barcelona Graduate School of Mathematics (BGSMath) and part of the Catalan research group 2017 SGR 01392.}

\maketitle


\vspace{0.3in}

\hspace{7cm} {\sl Dedicated to the memory of Pavol Brunovsk\'y}

\vspace{0.3in}

\begin{abstract} The purpose of this paper is to present an example of an Ordinary Differential Equation $x'=F(x)$ in the infinite-dimensional Hilbert space $\ell^2$ with $F$ being of class $\mathcal{C}^1$ in the Fr\'{e}chet sense, such that the origin is an asymptotically stable equilibrium point but the spectrum of the linearized operator $DF(0)$ intersects the half-plane $\Real(z)>0$. The possible existence or not of an example of this kind has been an open question until now, to our knowledge. An analogous example, but of a non-invertible map instead of a flow defined by an ODE  was recently constructed by the authors in \cite{8Ro-Sola}. The two examples use different techniques, but both are based on a classical example in Operator Theory due to S. Kakutani.

\end{abstract}

\noindent 

\noindent 

\section{Introduction}

In \cite{8Ro-Sola} we exhibited an example of a $\mathcal{C}^1$ map, in the Fréchet sense, in the Hilbert space $\ell^2$, namely $T:\ell^2\to\ell^2$,  such that $T(0)=0$,  the spectral radius of $DT(0)$ is larger than 1, but that the origin is exponentially asymptotically stable for the discrete nonlinear semi-dynamical system defined by the iterates $T^n$, for $n\ge 0$. By no means in that example $T$ is a diffeomorphism, since even the linear operator $DT(0)$ has no continuous inverse. In fact, $DT(0)$ is not onto. Therefore, in particular, the map $T(x)$ cannot be the time-one map of an Ordinary Differential Equation (ODE) in $\ell^2$.

That result was of recent interest because of the results and comments made by T. Gallay, B. Texier, and K. Zembrun in \cite{GTZ}, where they asked if a map like this could exist or not. They also pointed the difference with the case of Gateaux differentiability. After \cite{8Ro-Sola} was published, other interesting comments and results appeared in the paper \cite{GPP} by Á. Garab, M. Pituk, and C. Pötzsche. Perhaps the most relevant of these comments is that an example like that of \cite{8Ro-Sola} would not fall, and would not be possible, into the more stringent definition of exponential stability they use, that is slightly different from ours.

The goal of the present paper is to show that with a suitable adaptation of some of the ideas of  \cite{8Ro-Sola}, one can get an ODE $x'=F(x)$ in $\ell^2$ with a right hand side $F(x)$ of class $\mathcal{C}^1$, in the Fréchet sense, such that $F(0)=0$ and that the spectrum of $DF(0)$ intersects the half-plane $\Real (z)>0$, but the equilibrium point $x=0$ is asymptotically stable for the flow in the sense of Lyapunov. In fact, $x=0$ will be proven to be globally asymptotically stable.

We believe that the fact of being a flow defined by an ODE, and not just a map, is important for several reasons. First, because it somehow closes the implicit question appeared in the classical reference \cite{DK} by J.L. Daleckii, and M.G. Krein, about the discussion on how the spectrum of the linear part gives sufficient conditions for stability or instability. This discussion was originally presented, in that book, only for ODEs, and it was only later on, in the book \cite{Henry} by D. Henry, where the question was translated, also implicitly, to maps. 

In the case of ODEs $x'=F(x)$ in \cite{DK} and in the case of maps $x\mapsto T(x)$ in \cite{Henry} it was proved that when the nonlinearity is of class $\mathcal{C}^1$, and satisfies the extra regularity assumption $F(x)-DF(0)x=O(|x|^{1+p})$ in the first case, and $T(x)-DT(0)x=O(|x|^{1+p})$ in the second case, always with $p>0$, then the condition $\sigma(DF(0))\cap\{\Real(z)>0\}\ne\emptyset$, in the case of ODEs, or the condition $\sigma(DT(0))\cap\{|z|>1\}\ne\emptyset$ in the case of maps, imply instability without more additional assumptions. These extra regularity assumptions are satisfied, for example, when $F$ or $T$ are of class $\mathcal{C}^2$. Therefore, it is well-known that our examples could not be possible with nonlinearities of class $\mathcal{C}^2$ or of class $\mathcal{C}^1$ with the mentioned extra regularity hypothesis.

The other reason is that  we have seen in several contexts that the relation between the dynamics of a map $T$ in a neighborhood of the fixed point $x=0$ and that of $DT(0)$ is stronger when $DT(0)$ is invertible. This can be seen in our previous paper \cite{6Ro-Sola}, in relation to a different problem, namely that of smooth linearization. In fact, there are results that are known for diffeomorphisms that are still unknown, and maybe not true, for non-invertible smooth maps. Therefore, it is interesting to see that the time-one map of the ODE presented in the present paper is a diffeomorphism that exhibits the same properties as the map $T(x)$ of \cite{8Ro-Sola}.

We want to say that in the present paper we leave at least two open problems. The first is if the non-invertible map $T(x)$ of our previous paper  \cite{8Ro-Sola} could be or not the time-one map of a semi-flow, defined by a differential equation in $\ell^2$ whose linear part is an unbounded operator, 
and the second is if the asymptotic stability proved in the present paper is in fact or not an exponential stability, whatever definition we take, or, at least, if this exponential stability could be proved for a similar, but perhaps not exactly equal example of an ODE in $\ell^2$.
  
Our example in \cite{8Ro-Sola} was based on an example in Operator Theory, due to S. Kakutani (see also the books of Ch.E. Rickart \cite{Rickart} or P.R. Halmos \cite{Halmos}), that can be summarized as follows: there exists an operator $W_\varepsilon\in\mathcal{L}(\ell^2)$, depending on a certain sequence of real numbers $\varepsilon=(\varepsilon_m)_{m\ge 1}$,  with spectral radius $\rho(W_\varepsilon)>1$ and also a  sequence of operators $(L_m)_{m\ge 1}\subset\mathcal{L}(\ell^2)$  such that $\|L_m\|\to 0$ as $m\to\infty$,  and
\begin{equation}\label{nil}
(W_\varepsilon-L_m)^{p(m)}=0
\end{equation}
for a certain power $p(m)$. Consequently, $\rho(W_\varepsilon-L_m)=0$.

Then, the idea in \cite{8Ro-Sola}, was essentially (but without some details that at the end became relevant) to define the nonlinear map $T$ as $T(x)=W_\varepsilon x-L_mx$ for $r_1(m)\le |x|\le r_2(m)$ and a kind of suitable smooth interpolation for $r_2(m+1)\le |x|\le r_1(m)$, for suitable sequences $r_1(m)$ and $r_2(m)$ that approach zero as $m$ tends to $\infty$. Somehow, this is like a discontinuous {\em switching system} that has been smoothed out. In fact, our example started being a discontinuous system, and later on we realized that it could be smoothed out. To built this smooth interpolation without loosing the dynamical properties, we became aware that property \eqref{nil} could be extended to a stronger property, valid on a larger class $\Omega_m$ of operators $(W_\varepsilon-L'_m)\in\Omega_m$, namely that any arbitrary product of $p(m)$ factors 
\begin{equation}\label{setnil}
(W_\varepsilon-L'_{m,1})\cdot(W_\varepsilon-L'_{m,2})\cdots(W_\varepsilon-L'_{m,p(m)})=0,
\end{equation}
as long as they belong to that class, even if they don't commute. This can be described in words by saying that the set $\Omega_m\subset\mathcal{L}(\ell^2)$ is a nilpotent set of index at most $p(m)$. Then, a smooth interpolation between, say, $(W_\varepsilon-L_m)$ and $(W_\varepsilon-L_{m+1})$ can be build, because the families of operators $\Omega_m$ and $\Omega_{m+1}$ have in fact a non-empty intersection. This smooth interpolation is a delicate point, since we want the resulting function to be of class $\mathcal{C}^1$ even at $x=0$.

The idea in the present paper has to be different, because a nilpotent linear operator can never be obtained as a solution flow of an ODE. But the essential point in \eqref{setnil} is that the operators $(W_\varepsilon-L'_{m,i})$ have a norm perhaps larger than one, but the product of $p(m)$ of them has a very small norm, because it is zero. The idea in the present paper is to consider the ODE $x'=F(x)$ with a right-hand side $F(x)$ that is related to Kakutani's construction, and solve it by the method of Euler's polygons, that is $x_{i+1}=x_i+hF(x_i)$ and then prove, in analogy with \eqref{setnil}, that the composed map $(I+hF(\cdot))\circ (I+hF(\cdot))\circ\cdots (p(m) \text{ terms} )\cdots \circ (I+hF(\cdot))$
can be a contraction, for $m$ sufficiently large, and $h$ suitably small, even if each of the terms is not.

\section{Notation and the fundamental Lemma}

Let us start by explaining our notation. In the Hilbert space $\ell^2$, let $(e_n)_{n\ge 1}$ be the canonical Hilbert basis. A {\em weighted shift} $W\in\mathcal{L}(\ell^2)$, is an operator defined by $We_n=\alpha_ne_{n+1}$, and $(\alpha_n)$ is called its {\em weights sequence}.

The next is Kakutani's construction. Let $\varepsilon=(\varepsilon_m)_{m\ge 1}$, be a real sequence, strictly decreasing and tending to zero. For $n\ge 1$ let us define the integer $\kappa(n)\ge 0$ as $n=2^{\kappa(n)}(2\ell+1)$. Then, we define $\alpha_n=\varepsilon_{1+\kappa(n)}$.  Therefore, the member $\varepsilon_\mu$ of the $\varepsilon$ sequence appears at the places $\alpha_n$ for $n=2^{\mu-1}(2\ell+1)$ for $\ell=0,1\dots$. Therefore, it appears periodically with a period $2^\mu$ and appears for the first time at the position $n=2^{\mu-1}$.

The operator $W_\varepsilon$ is the weighted shift defined by the weights sequence $\alpha_n=\varepsilon_{1+\kappa(n)}$ obtained after the sequence $(\varepsilon_m)_{m\ge 1}$. The weighted shifts $L_k$ are defined by the same formula as $W_\varepsilon$ but with all the $\varepsilon_m$ replaced by $0$ except when $m=k$. Observe that $\|L_k\|\to 0$ as $k\to \infty$. The operator $W_\varepsilon-L_k$ is a nilpotent operator of index $2^k$ (exactly). Kakuntani's example consists of the observation that the spectral radius $\rho(W_\varepsilon)$ can be positive, while $\rho(W_\varepsilon-L_k)=0$ for all $k$.

But we observed in \cite{8Ro-Sola} that all weighted shift operators defined by a weights sequence $\alpha_n$ that vanishes at the same places, namely, when $n=2^{k-1}(2\ell+1)$ for $\ell=0,1\dots$ are also nilpotent of index (at most) $2^k$, an they constitute a class $\Omega_k\subset\mathcal{L}(\ell^2)$ with this important extra property: if $W_1,W_2,\dots W_{2^k}\in\Omega_k$, then $W_1\cdot W_2\cdots W_{2^k}=0$. To convince oneself of this property, the best is to see that for all the elements $e_n$ of the Hilbert basis one has $(W_1\cdot W_2\cdots W_{2^k})(e_n)=0$. 

For later reference, let us repeat, more synthetically, the definition of $\Omega_k$, ($k\ge 1$):
\begin{equation}\label{Omegak}
\Omega_k=\{W\in\mathcal{L}(\ell^2)\,|\, We_n=\alpha_ne_{n+1}\text{ and } 
\alpha_n=0 \text{ if } n=2^{k-1}(2j+1)\text{ for all } j=0,1\dots\}. 
\end{equation}

Let us see now that for a good choice of the sequence $(\varepsilon_m)$ the spectral radius $\rho(W_\varepsilon)$ of $W_\varepsilon$ turns out to be equal to 1. Observe that $\|W_\varepsilon^n\|=\alpha_1\alpha_2\cdots\alpha_n$ and that for integer $p\ge 1$ and $n=2^p-1$
$$\|W_\varepsilon^n\|=\|W_\varepsilon^{2^p-1}\|=\varepsilon_{1}^{2^{(p-1)}}\varepsilon_2^{2^{(p-2)}}\cdots\varepsilon_{p-2}^4\varepsilon_{p-1}^2\varepsilon_p^1=\prod_{q=1}^{p}\varepsilon_q^{2^{(p-q)}}.$$
Hence,
$$\log \left(\|W_\varepsilon^n\|^{\frac{1}{n}}\right)=\dfrac{1}{2^p-1}\sum_{q=1}^p2^{(p-q)}\log\varepsilon_q=\dfrac{2^p}{2^p-1}\sum_{q=1}^p \dfrac{\log\varepsilon_q}{2^q}.$$

Let us take now $(\varepsilon_m)$ as $\varepsilon_m=1/K^{m-2}$ for some $K>1$. Then, $\|W_\varepsilon\|=K>1$ and it is easy to calculate that $\rho(W_\varepsilon)=1$, with Gelfand's formula, since $\sum_{q=1}^\infty (q-1)2^{-q}=1$. Let us choose now $K=K_0>1$,  $\varepsilon_m=1/K_0^{m-2}$ and let us from now on fix the notation $W_\varepsilon$ to mean the fixed operator obtained with this choice $K=K_0$.

Observe that $\|W_\varepsilon\|, \|L_k\|,\|W_\varepsilon-L_k\|\le K_0$, for all $k$.

Now we are in the position to state and prove the next fundamental lemma:

\begin{lem}\label{fundamental}
Let us consider an ODE in $\ell^2$ of the form $x'=\bigl(-\gamma I+W(|x|)\bigr)x$, for some $\gamma>0$ where $W(s)$ is a Lipschitz function of $s$ in some compact interval $0<a\le s\le b<\infty$ into $\mathcal{L}(\ell^2)$. Suppose that for all $s$ in this interval, the operator $W(s)$ belongs to a single one of the classes $\Omega_k\subset\mathcal{L}(\ell^2)$ defined above, and suppose also that $\|W(s)\|\le K_1$, for all $s$. Suppose that $x(t)$, for $0\le t\le t_0$, satisfies $a\le |x(t)|\le b$ and is a solution of this ODE. Let $p_n(t)$ be the $n$-degree polynomial $p_n(t)=1+t+t^2/2+\cdots t^n/n!$. Then, we claim that
 \begin{equation}\label{lemma} 
 |x(t)|\le p_{(2^k-1)}(tK_1)e^{-\gamma t}|x(0)|,
\end{equation}
for all $0\le t\le t_0$.
\end{lem}

{\it Proof}: In an initial-value problem for an ODE of the form $x'=F(x)$ with  $x(0)=x_0$ the {\sl Euler-polygon} approximate solution with step $h$ is the piecewise-linear interpolation of some nodal values $(0,x_0),(h,x_1),(2h, x_2)\dots$ of the form $x^h(t)=x_n+(t/h-n)(x_{n+1}-x_n)$ for $nh\le t\le(n+1)h$, where the nodal values $x_n$ are defined by $x_{n+1}=x_n+hF(x_n)$. It is known that if $F$ is globally Lipschitz, the true solution exists for all time and the Euler polygons approach this true solution as $h\to 0$, uniformly on bounded intervals of the variable $t$ . Classical references for this property are the book of W. Hurewicz \cite{Hurewicz}, even if the proof is presented there only in two dimensions, and also the book of J. Dieudon\' ee \cite{Dieudonne}, statement (10.5.1.1) together with the remark (10.4.6).

Our function $F(x)=(-\gamma I+W(|x|))x$ is initially defined only for $a\le|x|\le b$ but can be extended to the whole of $\ell^2$ by defining $W(s)=W(a)$ for $0\le s<a$ and $W(s)=W(b)$ for $b<s$. With this extension $W(s)\in\Omega_k$ for all $s$ and also $F(x)$ is globally Lipschitz.

We are going to prove the inequality \eqref{lemma} for every solution of $x'=F(x)$, and we will remember that they will be solutions of the initial ODE provided that $a\le |x(t)|\le b$ for all $t\in[0,t_0]$, as it is said in the statement. 

Let us fix now the value of $t>0$. We divide the interval $[0,t]$ into $N$ parts, use the Euler polygon with $h=t/N$, and obtain \eqref{lemma} as a limit when $N\to \infty$ of a similar inequality obtained for the endpoint $x_N$.

If $F(x)=(-\gamma I+W(|x|))x$, then the nodes of the Euler polygon would be
$x_{i+1}=x_i+\frac{t}{N}(-\gamma x_i+W(|x_i|)x_i)$, for $i=0,1\dots (N-1)$. 

Therefore, \[x_N=\left(I+\frac{t}{N}(-\gamma I+W(|x_{N-1}|))\right)\cdot\left(I+\frac{t}{N}(-\gamma I+W(|x_{N-2}|))\right)\cdots\] \[ \left(I+\frac{t}{N}(-\gamma I+W(|x_{1}|))\right)\cdot\left(I+\frac{t}{N}(-\gamma I+W(|x_{0}|))\right)x_0.\]
The map $x_0\mapsto x_N$ is an approximation of the map $x(0)\mapsto x(t)$ that will improve as $N\to\infty$.

Written in other words
\[x_N=\left((1-\frac{t}{N}\gamma)I+\frac{t}{N}W(|x_{N-1}|)\right)\cdot\left((1-\frac{t}{N}\gamma)I+\frac{t}{N}W(|x_{N-2}|)\right)\cdots \] \[\left((1-\frac{t}{N}\gamma)I+\frac{t}{N}W(|x_{1}|)\right)\cdot\left((1-\frac{t}{N}\gamma)I+\frac{t}{N}W(|x_{0}|)\right)x_0.\]

To simplify, let us write $W_i$ instead of $W(|x_i|)$. Then
\[x_N=\left((1-\frac{t\gamma}{N})I+\frac{t}{N}W_{N-1}\right)\cdot\left((1-\frac{t\gamma}{N})I+\frac{t}{N}W_{N-2}\right)\cdots \left((1-\frac{t\gamma}{N})I+\frac{t}{N}W_0\right)x_0\]
and for $N\ge 2^k$ we have
\[x_N=\left((1-\frac{t\gamma}{N})I+\frac{t}{N}W_{N-1}\right)\cdot\left((1-\frac{t\gamma}{N})I+\frac{t}{N}W_{N-2}\right)\cdots \left((1-\frac{t\gamma}{N})I+\frac{t}{N}W_0\right)x_0\]
\[=\left ((1-\frac{t\gamma}{N})^NI+(1-\frac{t\gamma}{N})^{N-1}\frac{t}{N}\sum_{0\le i\le N-1} W_i+(1-\frac{\gamma}{N})^{N-2}\frac{t^2}{N^2}\sum_{0\le i_1<i_2\le N-1}W_{i_1}W_{i_2}+\right .\]
\[(1-\frac{\gamma}{N})^{N-3}\frac{t^3}{N^3}\sum_{0\le i_1<i_2<i_3\le N-1}W_{i_1}W_{i_2}W_{i_3}+\cdots +\]
\[\left .(1-\frac{\gamma}{N})^{N-2^k+1}\frac{t^{2^k-1}}{N^{2^k-1}}\sum_{0\le i_1<i_2<\dots <i_{(2^k -1)}\le N-1}W_{i_1}W_{i_2}\cdots W_{i_{(2^k -1)}}\right )x_0\]

\noindent with no more terms, since the products of more than $2^k-1$ factors of the $W_i$ become zero. Therefore,  if $N$ is sufficiently large to have $1-\frac{t\gamma}{N}>0$ (remember that $t$ is fixed), we obtain   
\begin{equation}\label{boundproduct1}
\begin{split}
\left\|\left((1-\frac{t\gamma}{N})I+\frac{t}{N}W_{N-1}\right)\cdot\left((1-\frac{t\gamma}{N})I+\frac{t}{N}W_{N-2}\right)\cdots \left((1-\frac{t\gamma}{N})I+\frac{t}{N}W_0\right)\right\|\le
\\ {N\choose 0}(1-\frac{t\gamma}{N})^N+ {N\choose 1}(1-\frac{t\gamma}{N})^{N-1}\frac{tK_1}{N}+\cdots + {N\choose 2^k-1}(1-\frac{t\gamma}{N})^{N-2^k+1}\frac{(tK_1)^{2^k-1}}{N^{2^k-1}}\le
\\ (1-\frac{t\gamma}{N})^{N-2^k+1}\left({N\choose 0}(1-\frac{t\gamma}{N})^{2^k-1}+ {N\choose 1}(1-\frac{t\gamma}{N})^{2^k-2}\frac{tK_1}{N}+\cdots + {N\choose 2^k-1}\frac{(tK_1)^{2^k-1}}{N^{2^k-1}}\right).
\end{split}
\end{equation}

It is a very elementary fact that ${N\choose i }\le\frac{N^i}{i!}$.  Also $0<1-\frac{t\gamma}{N}<1$ for $N$ sufficiently large. Therefore, the last expression in \eqref{boundproduct1} can be bounded by $(1-\frac{t\gamma}{N})^{N}(1-\frac{t\gamma}{N})^{1-2^k}  p_{(2^k-1)}(tK_1)$, where $p_{(2^k-1)}(z)$ is the polynomial of degree $2^k-1$ given by $p_{(2^k-1)}(z)=1+z+\frac{1}{2}z^2+\dots+\frac{1}{(2^k-1)!}z^{2^k-1}$. Taking now limits as $N\to\infty$ we get that $$|x(t)|\le p_{(2^k-1)}(tK_1)e^{-\gamma t}|x_0|,$$ for each $t>0$, as claimed. 
\cqd

\section{Construction of the example} 

Given two real numbers $a<b$ let $\varphi_1(a,b;s)$ be a real function of $s\in[a,b]$ of class $\mathcal{C}^1$ in $s$ such that $\varphi_1(a,b;a)=\partial_s\varphi_1(a,b;a)=0$, $\varphi_1(a,b;b)=1$, $\partial_s\varphi_1(a,b;b)=0$ and $0\le\varphi_1(a,b;s)\le 1$ for all $s\in[a,b]$. It is easy to see that  we can suppose also that $0\le\partial_s\varphi_1(a,b;s)\le 2/(b-a)$. Let us define also $\varphi_2=1-\varphi_1$, that reverses the values at the limits $a$ and $b$.

Let us define now  $\widetilde{L}_k(s)$, for $s\ge 0$ and $k\ge 3$
\begin{equation}\label{Lk}
  \widetilde{L}_k(s)=
  \begin{cases}
  0,\text{ if } s<r_{k+3}\\
  \varphi_1(r_{k+3},r_{k+2};s)L_k, \text{ if } s\in[r_{k+3},r_{k+2})\\
                                            L_k, \text{ if } s\in[r_{k+2},r_{k-1})\\
  \varphi_2(r_{k-1},r_{k-2};s)L_k, \text{ if } s\in[r_{k-1},r_{k-2})\\
   0,\text{ if } s\ge r_{k-2},

  \end{cases}
\end{equation}
for a suitable infinite sequence $(r_k)_{k\ge 1}$ strictly decreasing and tending to $0$, to be defined later. For the moment we only ask the condition 
\begin{equation}\label{geometric}
r_{k+1}\le r_k/2, \text{ for all } k\ge 1.
\end{equation}

See Fig. \ref{graph} for a sketch of the graph of $\widetilde{L}_k(s)$. It is clear that $\widetilde{L}_k:[0,\infty)\to\mathcal{L}(\ell^2)$ is a $\mathcal{C}^1$ function of compact support such that $W_\varepsilon-\widetilde{L}_k(s)\in\Omega_k$ for all $s\in [r_{k+2},r_{k-1}]$.

\begin{figure}[h] \label{graph}
\begin{tikzpicture}[xscale=10,yscale=1/3]
\draw[->,black, thin, domain=0:1.45] plot (\x,0);
\node at (1.50,-.05) {$s$};
\draw[black,thin] ( .15,0.4) -- ( .15,-0.4) node[below,black]{$r_{k+3}$};
\draw[black,thin] ( .34,0.4) -- ( .34,-0.4) node[below,black]{$r_{k+2}$};
\draw[black,thin] ( .55,0.4) -- ( .55,-0.4) node[below,black]{$r_{k+1}$};
\draw[black,thin] ( .78,0.4) -- ( .78,-0.4) node[below,black]{$r_{k}$};
\draw[black,thin] ( 1.03,0.4) -- ( 1.03,-0.4) node[below,black]{$r_{k-1}$};
\draw[black,thin] (1.30,0.4) -- (1.30,-0.4) node[below,black]{$r_{k-2}$};
\draw[black, thick, domain=0:0.15] plot (\x,.08+0);
\draw[black, thick, domain=0.15:0.34] plot (\x,{.08  +(5*((\x-0.15)/(.19))^2-10*(\x-0.15)^2*(\x-0.34)/(.19^3)});
\draw[black, thick, domain=0.34:1.03] plot (\x,.08+5);
\draw[black, thick, domain=1.03:1.30] plot (\x,{.08+5-(5*((\x-1.03)/(.27))^2-10*(\x-1.03)^2*(\x-1.30)/(.27^3)});
\draw[black, thick, domain=1.30:1.40] plot (\x,.08+0);
\node at (.690,7.0) {{\sl $\widetilde{L}_k(s)=L_k$}};
\node at (0.07,4.2) {{$\widetilde{L}_k(s)=0$}};
\node at (1.40,4.2) {{$\widetilde{L}_k(s)=0$}};
\draw[->,black,thin] (0.86,6.5) -- (0.93,5.5);
\draw[->,black,thin] (0.52,6.5) -- (0.45,5.5);
\draw[->,black,thin] (.07,3.2) -- (.07,1.2);
\draw[->,black,thin] (1.33,3.2) -- (1.33,1.2);
\end{tikzpicture}
\caption{Sketch of the graph of $\widetilde{L}_k(s)$}
\label{graph}
\end{figure}
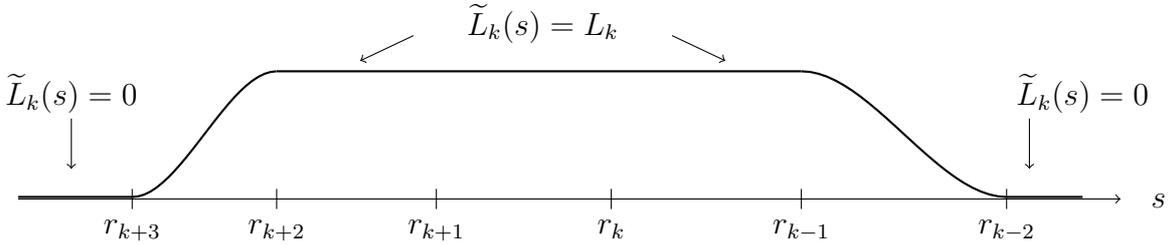

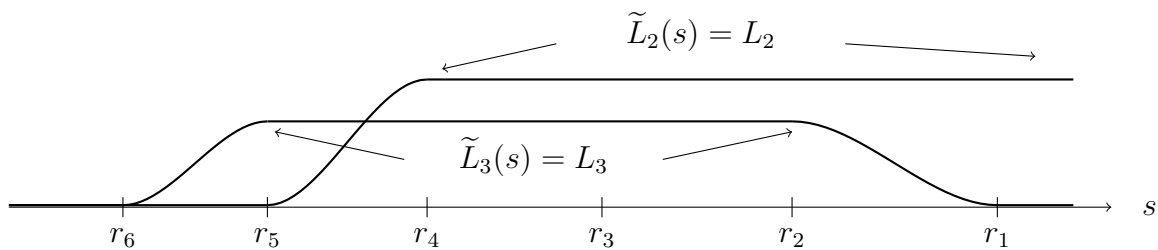
\begin{figure}[h]
\begin{tikzpicture}[xscale=10,yscale=1/3]
\draw[->,black, thin, domain=0:1.45] plot (\x,0);
\node at (1.50,-.05) {$s$};
\draw[black,thin] ( .15,0.4) -- ( .15,-0.4) node[below,black]{$r_6$};
\draw[black,thin] ( .34,0.4) -- ( .34,-0.4) node[below,black]{$r_5$};
\draw[black,thin] ( .55,0.4) -- ( .55,-0.4) node[below,black]{$r_4$};
\draw[black,thin] ( .78,0.4) -- ( .78,-0.4) node[below,black]{$r_3$};
\draw[black,thin] ( 1.03,0.4) -- ( 1.03,-0.4) node[below,black]{$r_2$};
\draw[black,thin] (1.30,0.4) -- (1.30,-0.4) node[below,black]{$r_1$};
\draw[black, thick, domain=0:0.15] plot (\x,.08+0);
\draw[black, thick, domain=0.15:0.34] plot (\x,{.08  +(5*((\x-0.15)/(.19))^2-10*(\x-0.15)^2*(\x-0.34)/(.19^3))/1.5});
\draw[black, thick, domain=0.34:1.03] plot (\x,.08+5/1.5);
\draw[black, thick, domain=1.03:1.30] plot (\x,{.08+((5-(5*((\x-1.03)/(.27))^2-10*(\x-1.03)^2*(\x-1.30)/(.27^3))))/1.5});
\draw[black, thick, domain=1.30:1.40] plot (\x,.08+0);
\node at (.690,7.0/3.5) {{\sl $\widetilde{L}_3(s)=L_3$}};
\draw[->,black,thin] (0.86,1.9) -- (1.03,3.0);
\draw[->,black,thin] (0.52,1.9) -- (0.35,3.0);
\draw[black, thick, domain=0:0.34] plot (\x,.08+0);
\draw[black, thick, domain=0.34:0.55] plot (\x,{.08  +(5*((\x-0.34)/(.21))^2-10*(\x-0.34)^2*(\x-0.55)/(.21^3)});
\draw[black, thick, domain=0.55:1.40] plot (\x,.08+5);
\node at (.91,7.0) {{\sl $\widetilde{L}_2(s)=L_2$}};
\draw[->,black,thin] (1.1,6.5) -- (1.35,6);
\draw[->,black,thin] (0.72,6.5) -- (0.57,5.5);
\end{tikzpicture}
\caption{Sketch of the graphs of $\widetilde{L}_2(s)$ (above) and $\widetilde{L}_3(s)$ (below)}
\label{graph2}
\end{figure}

Let us define  also $\widetilde{L}_2(s)$, for $s> 0$:
\begin{equation}\label{L2}
  \widetilde{L}_2(s)=
  \begin{cases}
  0,\text{ if } s<r_{5}\\
  \varphi_1(r_{5},r_{4};s)L_2, \text{ if } s\in[r_{5},r_{4})\\
                                            L_2, \text{ if } s\in[r_{4},\infty).

  \end{cases}
\end{equation}
See Fig. \ref{graph2}. The reader will find perhaps surprising that we do not define something like $\widetilde{L}_1(s)$ or even that the operator $L_1$ will not appear in our construction. This is made on purpose, though the reasons are not very deep: we don't want the sequence $(r_k)$ to start with $k=0$ and at the same time that we find convenient to keep a central place for the interval $[r_{k+1},r_k)$ in the definition of $\widetilde{L}_(s)$, as it appears in Fig. \ref{graph}.

Let us define now $\widetilde{L}(s)$, for $s\ge 0$, by  $\widetilde{L}(s)=\sum_{k\ge 2}\widetilde{L}_k(s)$. Observe that in this last infinite sum, a maximum of five terms are different from zero for each $s>0$. This implies that 
\begin{equation}\label{Lsbound}
\|\widetilde{L}(s)\|\le5K_0
\end{equation}
for all $s>0$. We also see that $\widetilde{L}(0)=0$, and the fact that $\|{L}_k\|\to 0$
as $k\to\infty$ implies that $\widetilde{L}(s)$ is continuous at $s=0$.

Given any $s_0\in(0,r_4)$ there exists one and only one $k\ge 4$ such that $s_0\in[r_{k+1},r_k)$. Then, for $s$ in some small neighborhood of $s_0$ the infinite sum $\sum_{\ell\ge 2}\widetilde{L}_\ell(s)$ consists of a maximum of only six nonzero terms, namely 
$\widetilde{L}_{k+3}(s)+\widetilde{L}_{k+2}(s)+\widetilde{L}_{k+1}(s)+\widetilde{L}_{k}(s)+\widetilde{L}_{k-1}(s)+\widetilde{L}_{k-2}(s)$. With obvious simplifications if $s_0\in[r_4,r_3)$ or $s\in[r_3,\infty)$.This implies that $\widetilde{L}(s)$ is of class $\mathcal{C}^1$ in all of $s>0$.

Now we want to bound the function $s\widetilde{L}'(s)$. We first recall from the definition of $\varphi_1$ that $0<\partial_s\varphi_1(r_{k+3},r_{k+2},s))\le2/(r_{k+2}-r_{k+3})$. Then, if $s\in[r_{k+3},r_{k+2})$ we have that $\|s\widetilde{L}_k'(s)\|\le r_{k+2}2/(r_{k+2}-r_{k+3})\|L_k\|\le 4\|L_k\|$ because of property \eqref{geometric}, and this last quantity tends to $0$ as $k\to\infty$. The case $s\in[r_{k-1},r_{k-2})$ is similar, namely $\|s\widetilde{L}_k'(s)\|\le 4\|L_k\|$, and in the other cases $\widetilde{L}'_k(s)=0$. Summarizing, $\|s\widetilde{L}_k'(s)\|\le 4\|L_k\|$ for all $s>0$. This is for $k\ge 3$, but the bound for $\|s\widetilde{L}_2'(s)\|$ can be taken to be the same.

As we said above, given any $s_0\in(0,r_4)0$ there exists a $k$ such that $s_0\in[r_{k+1},r_k)$, and $\sum_{\ell\ge 2}\widetilde{L}_\ell(s)=\sum_{k-2\le\ell\le k+3}\widetilde{L}_\ell(s)$ for $s$ in some small neighbourhood of $s_0$. Therefore, in that small neighborhood of $s_0$, $\widetilde{L}'(s)=\sum_{k-2\le\ell\le k+3}\widetilde{L}'_\ell(s)$, and, dropping all the derivatives that are zero, $\widetilde{L}'(s_0)=\widetilde{L}'_{k+2}(s_0)+\widetilde{L}'_{k-2}(s_0)$. Then, because of the bound of the previous paragraph, we obtain that
\begin{equation}\label{sL'}
\|s_0\widetilde{L}'(s_0)\|\le 4\|L_{k+2}\|+4\|L_{k-2}\|\text{ when } r_{k+1}\le s_0<r_k.
\end{equation}

Then, since $k\to\infty$ when $s_0\to 0$, the previous bound implies that $s\widetilde{L}'(s)\to 0$ as $s\to 0$ and we conclude that the map $x\mapsto \widetilde{L}(|x|)x$ is of class $\mathcal{C}^1(\ell^2,\ell^2)$ because of the following lemma: 

\begin{lem}
	Let $H$ be a Hilbert space and $\widetilde{L}:[0,\infty)\to \mathcal{L}(H)$ be a continuous function such that  $\widetilde{L}(0)=0$. Suppose that $\widetilde{L}$ is of class $\mathcal{C}^1$ when restricted to $(0,\infty)$ and that $\ s\widetilde{L}'(s)\ $ tends to zero as $s$ tends to zero.
	Let $N:H\to H$ be defined by $N(x)=\widetilde{L}(|x|)x$. We claim that $N$ is of class $\mathcal{C}^1$, in the Fréchet sense, and $DN(0)=0$.
\end{lem}
{\sl Proof:} Let us start by showing that $DN(0)=0$ in the Fréchet sense:
$$\left|\dfrac{N(0+y)}{|y|}\right|=\left|\dfrac{\widetilde{L}(|y|)y}{|y|}\right|\le \|\widetilde{L}(|y|)\|,$$
and this last expression tends to $0$ as $|y|$ tends to $0$.

Let us now observe that the function $H\ni x\mapsto |x|$ is of class $\mathcal{C}^\infty$, in the Fréchet sense, in $H\setminus\{0\}$ and its first derivative at a point $x_0\ne 0$ is the linear map $y\mapsto D|x_0|[y]=\frac{\langle x_0,y\rangle}{|x_0|}$. The only thing it remains to prove is that $DN(x)$ approaches $0$ in $\mathcal{L}(H)$ as $x$ tends to $0$ in $H$:
$$\sup_{|y|=1}|DN(x)[y]|=\sup_{|y|=1}\left|\widetilde{L}'(|x|)\frac{\langle x,y\rangle}{|x|}x+\widetilde{L}(|x|)y\right|\le \|\widetilde{L}'(|x|)\||x|+\|\widetilde{L}(x)\|,$$
and these last quantities tend to $0$ as $|x|$ tends to $0$.\cqd

An important consequence of \eqref{sL'} is that
\begin{equation}\label{sL'2}
\|s\widetilde{L}'(s)\|\le 8K_0
\end{equation}
for all $s>0$, independently of the interval to which $s$ belongs. Strictly speaking, the bound was derived in \eqref{sL'} only for $s<r_4$, but can be extended easily to $s>r_4$.

Concerning the properties related to the sets $\Omega_k$, defined in \eqref{Omegak}, from the definition of the $\widetilde{L}_k(s)$ , for $k\ge 3$ we see that $W_\varepsilon-\widetilde{L}_k(s)=W_\varepsilon-L_k\in\Omega_k$ when $r_{k+2}\le s\le r_{k-1}$ and also that $W_\varepsilon-\widetilde{L}(s)\in\Omega_k$ when $r_{k+2}\le s\le r_{k-1}$. This, in particular, implies that if $r_{k+1}\le s\le r_{k}$ then $W_\varepsilon-\widetilde{L}(s)\in\Omega_{k+1}\cap
\Omega_{k}\cap\Omega_{k-1}$.

Let us consider now the $\mathcal{C}^1$ nonlinear map $F:\ell^2\to\ell^2$ defined by $F(x)=\bigl(-\gamma I+W_\varepsilon-\widetilde{L}(|x|)\bigr)x$, where we choose $\gamma$ with $0<\gamma<1$, and we fix it from now on, and consider also the ODE $x'=F(x)$ in $\ell^2$. Let us estimate its Lipschitz constant. 

Let $G(x)= \tilde L(|x|)x$. When $x,y\in\ell^2$ are not collinear, then the map $\theta\mapsto |\theta x+(1-\theta)y|$ is smooth and we can apply the chain rule in the next argument: 
$$|G(x)-G(y)|=\left|\int_0^1\frac{d}{d\theta}G(\theta x +(1-\theta )y)\ d \theta\right| 
$$
$$=\left|\int_0^1\widetilde{L}'(|\theta x +(1-\theta )y|)\frac{1}{2|\theta x +(1-\theta )y|}2\langle\theta x +(1-\theta )y,[x-y]\rangle(\theta x +(1-\theta )y)\ d\theta\right.$$
$$\left.+\int_0^1\widetilde L(|\theta x +(1-\theta )y|)[x-y]\ d\theta\right|\leq
\left(\int_0^1\| \widetilde{L}'(|\theta x +(1-\theta )y|)\| |\theta x +(1-\theta )y|\ d\theta\right) |x-y| 
$$
$$+\left(\int_0^1\widetilde{L}(|\theta x +(1-\theta )y|)\ d\theta\right) |x-y|.$$
Observe that in the first integral the subintegral function has the form $\|s\widetilde{L}'(s)\|$ and because of $\eqref{sL'2} $ this can be bounded by $8K_0$, and that in the second integral the subintegral function $\|\widetilde{L}(s)\|$ can be bounded by $5K_0$, because of \eqref{Lsbound}. Continuing then the previous inequality we obtain that 
$$|G(x)-G(y)|\le 12K_0|x-y|.$$
This last inequality has been proved when $x,y\in\ell^2$ are not collinear, but can be extended to the whole $(x,y)\in \ell^2\times\ell^2$ by continuity.

In conclusion, $F(x)$ is globally Lipschitz, that implies existence and uniqueness of solutions defined for all $t\in\R$.

We observe that $DF(0)=-\gamma I+W_\varepsilon$ and since the spectrum of $W_\varepsilon$ includes the set $|\lambda |=1$ we see that the spectrum of $DF(0)$ intersects the half-plane $\Real(z)>0$. This implies that the origin is unstable for the linearized ODE $x'=DF(0)x$.

Let us prove now that $x=0$ is asymptotically stable for the flow defined by $x'=F(x)$.
The proof will be carried on in two steps, and the proof of each step will be done by contradiction. After that, we will prove that the asymptotic stability is global.

In the first step we will prove that for all $k\ge 3$, if $r_{k+1}\le |x(0)|<r_k$ then $|x(t)|<r_{k-1}$ for all $t\ge 0$. To prove that, we will have to be precise about the definition of the sequence $(r_k)$, that, until now, was only restricted by \eqref{geometric}. At the end of this first step we conclude that $x=0$ is stable in the sense of Lyapunov.

In the second step we will prove, for all $k\ge 3$, that if $r_{k+1}\le|x(0)|<r_k$ then there exists a $t_1>0$ such that $|x(t_1)|< r_{k+2}$. This, together with the first step, proves asymptotic stability, namely proves that $|x(t)|\to 0$ as $t\to\infty$. Let us say why: because of the first step one will then have $|x(t)|<r_{k+1}$ for $t\ge t_1$, and for the second step again, but now applied for initial conditions in the interval $[r_{k+2},r_{k+1})$, a $t_2>t_1$ will exist such that $x(t_2)< r_{k+3}$, and then $|x(t)|<r_{k+2}$ for all $t\ge t_2$. And so on.

To define the sequence $(r_k)$ we first recall the $n$-degree polynomial $p_n(t)=1+t+t^2/2+\cdots t^n/n!$ we have considered in the Lemma \ref{fundamental} above. Let us consider the function $q_n(t)=p_n(K_0t)e^{-\gamma t}$. This function, for $t\ge 0$ starts with $q_n(0)=1$, initially increases, since $q_n'(0)=K_0-\gamma>0$, and finally tends to $0$ as $t\to\infty$. Let us call $Q_n=\max_{t\ge 0}q_n(t)$, define $r_1$  arbitrarily, and define $r_{k}=\min\{r_{k-1}/2,r_{k-1}/(2Q_{2^k-1})\}$, that obviously satisfies \eqref{geometric}.

Let us now prove the statement of the first step. Suppose it does not hold, and for some $s_1>0$ one has $|x(s_1)|=r_{k-1}$ and $|x(t)|<r_{k-1}$ for all $t<s_1$. Then, there should exist an $s_0\in(0,s_1)$ such that $|x(s_0)|=r_k$ and $|x(t)|\in(r_k,r_{k-1})$ for $s_0<t<s_1$. Along this piece of solution the term $\widetilde{L}_k(s)$ would be active in the definition of $\widetilde{L}(|x(t)|)$, and, more precisely, $\widetilde{L}_k(|x(t)|)=L_k$, and $W_\varepsilon-\widetilde{L}(|x(t)|)$ would belong to $\Omega_k$. Then, Lemma \ref{fundamental} applied to this piece of solution would give
$|x(t)|\le p_{2^k-1}(K_0(t-s_0))e^{-\gamma(t-s_0)}|x(s_0)|$ and in particular $r_{k-1}=|x(s_1)|\le r_kQ_{2^k-1}$, that contradicts the inequality $r_{k-1}\ge 2r_kQ_{2^k-1}$ that is deduced from the definition of $r_k$ above.

Let us finally prove the second step. We proceed also by contradiction. Suppose that 
$r_{k+1}\le|x(0)|<r_k$. Because of the first step, we will then have  $|x(t)|<r_{k-1}$ for all $t\ge 0$. Suppose now that the statement of the second step does not hold. This means that for all $t\ge 0$ we will have $r_{k+2}\le |x(t)| < r_{k-1}$. But for this range, $r_{k+2}\le s < r_{k-1}$, one has that $W_\varepsilon-\widetilde{L}(s)\in \Omega_k$. Therefore, we can apply Lemma \ref{fundamental} and obtain that $|x(t)|\le p_{2^k-1}(K_0t)e^{-\gamma t}|x(0)|$ for all $t\ge 0$, which contradicts, for large values of $t$, the inequality $r_{k+2}\le |x(t)| $ written above.

The previous reasoning shows that if $|x(0)|<r_3$ then $\lim_{t\to\infty}x(t)=0$, because in that case $r_{k+1}\le|x(0)|<r_k$ for some $k\ge 3$. Let us see that the same holds if
$|x(0)|\ge r_3$. As long as $|x(t)|\ge r_3$, one will have, because of \eqref{L2}, $\widetilde{L}_2(|x(t)|)=L_2$ and $W_\varepsilon-\widetilde{L}(|x(t)|)\in\Omega_2$. See Fig. \ref{graph2}, to clarify that.

Therefore, we can apply Lemma \ref{fundamental} as long as $|x(t)|\ge r_3$ and obtain that $|x(t)|\le p_{2^2-1}(K_0t)e^{-\gamma t}|x(0)|$. This implies that the interval of values of $t$ for which $|x(t)|\ge r_3$ cannot be arbitrarily long, and that after some time $t_1$ one will have $|x(t_1)|< r_3$. From this time $t_1$ onwards one can apply the arguments of the previous paragraphs, and obtain $\lim_{t\to\infty}x(t)=0$, as we were willing to prove.

\end{document}